\documentclass[11pt,a4paper]{amsart}
\usepackage{times} % assumes new font selection scheme installed
\usepackage{amsmath} % assumes amsmath package installed
\usepackage{amssymb}  % assumes amsmath package installed
\usepackage{amsthm}
\usepackage{latexsym}
\usepackage{amsfonts}
\usepackage{xcolor}
\usepackage{ulem}
\usepackage[pdftex,pagebackref=false]{hyperref} % with basic options
\hypersetup{
    plainpages=false,       % needed if Roman numbers in frontpages
    unicode=false,          % non-Latin characters in Acrobat’s bookmarks
    pdftoolbar=true,        % show Acrobat’s toolbar?
    pdfmenubar=true,        % show Acrobat’s menu?
    pdffitwindow=false,     % window fit to page when opened
    pdfstartview={FitH},    % fits the width of the page to the window
    pdfnewwindow=true,      % links in new window
    colorlinks=true,        % false: boxed links; true: colored links
    linkcolor=blue,         % color of internal links
    citecolor=red,        % color of links to bibliography
    filecolor=magenta,      % color of file links
    urlcolor=cyan           % color of external links
}
\usepackage{graphicx}
\newtheorem{thm}{Theorem}[section]
\newtheorem{cor}[thm]{Corollary}
\newtheorem{lem}[thm]{Lemma}
\newtheorem{prop}[thm]{Proposition}
\newtheorem{defn}[thm]{Definition}
\newtheorem{rem}[thm]{Remark}
\newtheorem{example}[thm]{Example}

\numberwithin{equation}{section}
\allowdisplaybreaks
\newcommand{\R}{\ensuremath{\mathbb R}}    % Reelle Zahlen
\newcommand{\C}{\ensuremath{\mathbb C}}    % Komplexe Zahlen
\newcommand{\N}{\ensuremath{\mathbb N}}    % Nat"urliche Zahlen
\newcommand{\<}{\langle}
\renewcommand{\>}{\rangle}

\newcommand{\calH}{\mathcal H}

\newcommand{\calL}{\mathcal L}         
         
\newcommand{\calN}{\mathcal N}         
         
\newcommand{\calP}{\mathcal P}

\newcommand{\calS}{\mathcal S}         
\newcommand{\calT}{\mathcal T}         
\newcommand{\calU}{\mathcal U}

\newcommand{\calY}{\mathcal Y}

\newcommand{\la}{\lambda}

\newcommand{\eps}{\varepsilon}

\newcommand{\bmat}[4]
{
	\begin{bmatrix}
		#1 & #2\\
		#3 & #4
	\end{bmatrix}
}
\newcommand{\bvec}[2]
{
	\begin{bmatrix}
		#1\\
		#2
	\end{bmatrix}
}

\renewcommand{\Re}{\operatorname{Re}}

\newcommand{\ran}{\operatorname{ran}}

\newcommand{\sra}{\rightarrow}

\newcommand{\Sra}{\Rightarrow}

\newcommand{\ol}{\overline}

\newcommand{\wt}{\widetilde}

\newcommand{\ptl}{\partial}

\newcommand\resout{\bgroup\markoverwith
{\textcolor{magenta}{\rule[.5ex]{2pt}{0.4pt}}}\ULon}

\begin{document}

\title{Necessary conditions for turnpike property for generalized linear-quadratic problems}
\author[R.\ Guglielmi]{Roberto Guglielmi}
\author[Z.\ Li]{Zhuqing Li}

\maketitle
\begin{abstract}
In this paper, we develop several necessary conditions of turnpike property for generalizaid linear-quadratic (LQ) optimal control problem in infinite dimensional setting. The term 'generalized' here means that both quadratic and linear terms are considered in the running cost.
The turnpike property reflects the fact that over a sufficiently large time horizon, the optimal trajectories and optimal controls stay for most of the time close to a steady state of the system. We show that the turnpike property is strongly connected to certain system theoretical properties of the control system. We provide suitable conditions to characterize the turnpike property
in terms of the detectability and stabilizability of the system. Subsequently, we show the equivalence
between the exponential turnpike property for generalized LQ and LQ optimal control problems.

\end{abstract}

\section{Introduction}
The turnpike phenomena have been observed and extensively studied in the context of mathematical economics in the early works by Von Neumann~\cite{VonNeumann} in 1945, which highlighted the trend
of optimally growing economies to approach balanced equilibrium path over a sufficiently long
time horizon. In 1958, Dorfman, Samuelson and Solow coined the term 'turnpike property'~\cite{Samuelson}. Over the last decade, various turnpike properties have been defined and extensively studied by
the mathematical community in the context of optimal control, also in connection with the stability of Model Predictive Control (MPC) schemes~\cite{LGMPC} and the qualitative properties of the control systems~\cite{Zua}. In general, the turnpike property describes the long-time
behavior of the optimally controlled system, where the optimal trajectories and controls over a
sufficiently large time horizon stay close to a particular state or trajectory of the system for most
of the time. Turnpike phenomena have received increasing interest in the control field due to the
structural insights they provide into the structure of optimal solutions. See, e.g., McKenzie~\cite{McKenzie}. For instance, on one hand, the turnpike property helps to extend the stability
results~\cite{AngeliD} obtained by constructing an appropriate Lyapunov function using strict dissipativity for the MPC closed
loop to larger classes of MPC schemes; on the other hand, the occurrence of turnpike property is
known to be closely linked to some structural-theoretical properties of the system like detectability and stabilizability. The turnpike property can also be used to synthesize long-term optimal
trajectories~\cite{Anderson,GTZ}. The monographs~\cite{ZAS06,ZAS14,ZAS15} present a complete overview on turnpike properties in various optimal control and variational problems.

This extensive literature on the subject has given rise to different notions of the turnpike
property, to a large extent related to the methodology adopted in different works. In particular, the measure turnpike and the exponential turnpike are two notions of particular importance.
The notion of exponential turnpike property naturally appears when exploiting the hyperbolicity feature around the turnpike state of the optimality system resulting from the Pontryagin’s
maximum principle~\cite{Zua,TreZua,TreZuaFIN}. On the other hand, the measure turnpike property
emerges when investigating the connection between turnpike and dissipativity properties of the
control system~\cite{GRUNE201645}. In most cases, the prescribed trajectory and control which is approached by the optimal pair is the minimizer of the corresponding optimal steady state problem, i.e., the optimal steady state of the system. However, turnpike phenomena have been observed towards
different target states, for example towards suitable periodic orbits ~\cite{SAMUELSONNONLINEAR,TreZua,Larsperiodic}.

Various sufficient conditions for turnpike properties have been derived in both the linear and
nonlinear setting. We refer to, e.g.,~\cite{LarsDiscrete,LarsDiscrete2} for discrete-time systems, \cite{Zua,TreZuaFIN} for finite dimensional continuous-time systems, \cite{Carlson} for delay differential equations, and \cite{arXiv,CARLSONDinfin,TreZua,Zua} for infinite
dimensional systems with control inside the domain. In comparison, only a few papers investigate necessary conditions for turnpike properties, mostly in the finite dimensional case. Indeed,
necessary conditions for measure turnpike property have been thoroughly studied in~\cite{RGDis,RGCon} for finite dimensional generalized LQ optimal control problems in the discrete-time and continuous-time setting, respectively. The combination of these findings with the sufficient conditions for
the turnpike property eventually yields a complete characterization of the turnpike property in
finite dimensions, mostly in terms of stabilizability and detectability of the control system.

In this work we aim at investigating the relation between turnpike property and certain notions of stabilizability and detectability for generalized linear quadratic optimal control problems
in the infinite dimensional setting. Specifically, we derive several independent necessary conditions for the
turnpike property: Statements (a) and (c) in Theorem~\ref{necess1} show that the exponential stabilizability of the control system and the exponential stability on the unobservable subspace are needed to guarantee the turnpike property, which provides a generalization of the finite-dimensional results in~\cite{RGCon} to the infinite dimensional setting. Furthermore, for the class of problems under consideration, Theorem~\ref{necess1} (b) resolves the question of the characterization of the turnpike reference by showing that the only possible turnpike reference is the optimal steady state, as defined in~\eqref{OSSP}. Moreover, we show that the necessary conditions in Theorem~\ref{necess1} become sufficient conditions of the turnpike property in the case of point spectrum (Theorem~\ref{suffnece1}) or finite dimensional state space (Corollary~\ref{suffnece2}). In particular, this provides a generalization of the finite dimensional counterpart in~\cite{RGCon}. Finally, in Theorem~\ref{suffnece3} we establish the equivalence between the exponential turnpike property for generalized LQ and
LQ optimal control problems.

Our paper is organized as follows. In Section \ref{Setting}, we introduce the problem setting of generalized LQ optimal control problem and define two notions of the turnpike property. In Section \ref{Results}, we present our main results and give an illustrative example. Section~\ref{Proofs} collects the proofs of
the results of the paper.

\section{Mathematical setting}\label{Setting}

\noindent
{\bf Notations}\\
We denote by $\calH,\calU,\calY$ the state, input and output spaces, respectively. We assume they are all complex Hilbert spaces.
We use $\|\cdot\|$ (resp. $\<\cdot,\cdot\>$) to denote the norm (resp. inner product) on all these spaces. For the sake of brevity, we denote by $\|\cdot\|_{L^2}$ (resp. $\<\cdot,\cdot\>_{L^2}$) the norm (resp. inner product) corresponding to the space $L^2(0,T;\calH)$, $L^2(0,T;\calU)$ and $L^2(0,T;\calY)$.

We are concerned with the following generalized LQ optimal control problem over the time interval $[0,T]$, with corresponding time horizon $T>0$:

\noindent\textbf{The generalized LQ optimal control problem} \noindent${(GLQ)}_T$. Find the optimal pair 
$$
(x^*(\cdot),u^*(\cdot))\in L^2(0,T;\calH)\times L^2(0,T;\calU)
$$
which minimizes the cost functional
\begin{equation}\label{OCP}
J_T(x_0,u):=\int_0^T\ell(x(t),u(t))dt,
\end{equation}
where the running cost $\ell : \mathcal{H}\times \mathcal{U}\to \R$ is defined by
\begin{equation}\label{cost}
\ell(x,u) := \|Cx\|^2 + \|Ku\|^2 + 2\Re\< z,x\> + 2\Re\< v,u\>,
\end{equation}
and the pair $(x,u)$ is constrained to the dynamical system
\begin{equation}\label{sy}
\addtolength{\jot}{2pt}
\left\{
\begin{alignedat}{2}
&\dot x(t)=Ax(t)+Bu(t),\\
&x(0)=x_0\in\calH.
\end{alignedat}
\right.
\end{equation}
Here, we assume that $C\in \calL(\calH,\calY),\,K\in \calL(\calU),\,z\in \mathcal{H}$ and $v\in \mathcal{U}$. The operator $K$ is further assumed to be coercive, i.e., there exists a constant $m>0$ such that $\<K^*Ku,u\>\geq m\|u\|^2$ for any $u\in \calU$.

Regarding the dynamical system~\eqref{sy}, we assume that $A:D(A)\to\calH$ is the generator of a strongly continuous semigroup $\calT = \left(\calT_t\right)_{t\ge 0}$ densely defined on $\calH$ and $B\in \calL(\calU, \calH)$. In particular, if $z=0$ and $v=0$, we shall refer to this problem as the LQ optimal control problem, which is abbreviated to $(LQ)_T$.

It is well-known that the optimal pair for problem $(GLQ)_T$ exists and is unique. See, e.g., \cite[Section 3.1]{thesis} for a proof. In the remainder of this paper, we denote by $x^*_T(\cdot,x_0)$ and $u^*_T(\cdot,x_0)$, or simply by $x^*$ and $u^*$ when $x_0$ and $T$ are clear from the context, the optimal trajectory and optimal control of problem $(GLQ)_T$ (or $(LQ)_T$) corresponding to initial condition $x_0\in \calH$ and time horizon $T>0$. 
\begin{rem}
A common choice for the running cost $\ell$ is the form of tracking problems
\begin{align*}
\ell(x,u):&= \|Cx-y\|^2+\|Ku-f\|^2\\
&=\|Cx\|^2+\|Ku\|^2-2\Re\<x,C^*y\>-2\Re\<u,K^*f\>+\|y\|^2+\|f\|^2\, .
\end{align*}
In this case, the running cost decreases to $0$ as the observation pair $(Cx,Ku)$ converges to the tracking reference $(y,f)\in \calY\times\calU$. Comparing this expression with the running cost~\eqref{cost} of a generalized LQ optimal control problem, the main difference is given by the fact that, in the running cost of tracking type, the pair $(z,v)$ is given by $(C^*y,K^*f)\in \ran C^*\times \ran K^*$; therefore, running cost of tracking type for a point reference can be seen as a special case of the running cost~\eqref{cost} associated with a generalized LQ optimal control problem. In the general case, we do not need to implicitly require the condition $(z,v)\in\ran C^*\times \ran K^*$, which allows us to consider a class of more general cost functionals. 
\end{rem}

\noindent{\bf Some notions of turnpike}\\
Following~\cite{Faulwasser,Zua}, we introduce the notions of the measure turnpike property and the exponential turnpike property at some steady state $(x_e,u_e)$. We recall that a steady state of system~\eqref{sy} is a pair $(x_e,u_e)\in D(A)\times\mathcal{U}$ such that $Ax_e+Bu_e=0$.

\begin{defn}\label{turnpikedef}
We say that the optimal control problem $(GLQ)_T$ satisfies the measure turnpike property at some steady state $(x_e,u_e)$ if, for any bounded neighborhood $\calN$ with respect to $\calH$-norm of $x_e$ and $\eps>0$, there exists a constant $M_{\calN,\eps}>0$ such that for all $x_0\in \calN$ and time horizon \(T>0\), the optimal trajectory $x_T^*(\cdot,x_0)$ and optimal control $u_T^*(\cdot,x_0)$ of problem $(GLQ)_T$ satisfy that
$$
\mu\left\{t\in[0,T]\;\big|\; \|x^*_T(t,x_0)-x_e\|+\|u^*_T(t,x_0)-u_e\|>\eps\right\}\leq M_{\calN,\eps},
$$
where $\mu$ denotes the Lebesgue measure on $\R$.

We say that the optimal control problem $(GLQ)_T$ satisfies the exponential turnpike property at some steady state $(x_e,u_e)$ if, for any bounded neighborhood $\calN$ of $x_e$, there exists some positive constants $M_{\calN}$ and $k$ such that for all $x_0\in \calN$ and time horizon \(T>0\), the optimal trajectory $x_T^*(\cdot,x_0)$ and optimal control $u_T^*(\cdot,x_0)$ of problem $(GLQ)_T$ satisfy
$$
\|x^*_T(t,x_0)-x_e\|+\|u^*_T(t,x_0)-u_e\|\leq M_{\calN}(e^{-kt}+e^{-k(T-t)}),\;\forall t\in[0,T].
$$
\end{defn}
\begin{rem}
It is clear that exponential turnpike property is stronger than the measure turn-
pike property. So, any sufficient condition for exponential turnpike property is automatically a
sufficient condition for measure turnpike property, and conversely, any necessary condition for
measure turnpike property is also a necessary condition for exponential turnpike property.
\end{rem}
A crucial element in our analysis of the turnpike property is the optimal steady state corresponding to the running cost $\ell$. The optimal steady state problem is defined as follows:
\begin{align}\label{OSSP}
\inf_{x \in D(A),\,u\in\calU} \ell(x,u)\quad\text{s.t. }Ax+Bu = 0.
\end{align}
If $(x_e,u_e)$ is a minimizer of optimal control problem~\eqref{OSSP}, we say $(x_e,u_e)$ is an optimal steady state.

We will also need several structural-theoretical properties of the control system under consideration, which we introduce in the following.

\begin{defn}
We say the pair $(A,B)$ is exponentially stabilizable if there exists some $F\in\calL(\calH,\calU)$ such that $A+BF$ generates an exponentially stable semigroup. We say the pair $(A,C)$ is exponentially detectable if the pair $(A^*,C^*)$ is exponentially stabilizable.
\end{defn}

Finally, let the operator $[A\;B]:D(A)\times \calU\to \calH$ be defined by
$$
[A\;B]\bvec{x}{u}:=Ax+Bu.
$$
Notice that, since $A$ and $B$ are both closed operators, $\ker[A\;B]$ is closed in $\calH\times\calU$.

We denote by $\bvec{A^*}{B^*}$ the adjoint operator of $[A\;B]$. It's easy to check that the domain of $\bvec{A^*}{B^*}$ is $D(A^*)$ and we have
\begin{align*}
\bvec{A^*}{B^*}w=\bvec{A^*w}{B^*w}{\in \calH\times\calU}\;,\quad \forall w\in D(A^*).
\end{align*}
\section{Main results}\label{Results}
In this section, we derive several necessary conditions for the turnpike property. Moreover,
we show that the turnpike property can be completely characterized by exponential stabilizability and detectability for the point spectrum case, and the exponential turnpike property for
generalized LQ optimal control problem and LQ optimal control problem are equivalent. The
proof of our results will be presented in the next section.

\begin{defn}
Given $C\in\calL(\calH,\calY)$, the unobservable subspace $U^{\infty}$ of a strongly continuous semigroup $\calT$ on $\calH$ is defined as
$$
U^{\infty}:=\{x\in\calH\;|\;C\calT_tx=0\;\;\text{for all}\;\; t\in [0,\infty)\}.
$$
\end{defn}
\begin{rem}
It is easy to check that $U^{\infty}$ is a closed subspace of $\calH$ and is invariant under $\calT$, so the restriction of $\calT$ to $U^{\infty}$ is a strongly continuous semigroup on $U^{\infty}$, and the corresponding generator is just $A|_{D(A)\cap U^{\infty}}$. See, e.g., ~\cite[Proposition 2.4.3]{TusWei}.
\end{rem}

Our first result provides several necessary conditions for the turnpike property in terms of the
turnpike reference, stabilizability and detecability of the system.

\begin{thm}\label{necess1}
If the problem $(GLQ)_T$ satisfies the measure or exponential turnpike property at some steady state $(x_e,u_e)$, then following statements hold:
\begin{enumerate}
\item[(a)] $\calT$ is exponentially stable on $U^{\infty}$.
\item[(b)] $(x_e,u_e)$ is the unique optimal steady state of problem \eqref{OSSP}.
\item[(c)] The pair $(A,B)$ is exponentially stabilizable.
\end{enumerate}
\end{thm}

Now assume $0\leq \tau<T$. Let $U_T(\cdot,\tau):[\tau,T]\to\calL(\calH)$ denote the evolution operator of the following problem: 
\begin{equation}\label{evoopttra}
\dot x(t)=(A-B(K^*K)^{-1}B^*P(T-t))x(t),\quad x(\tau)=x_0\in \calH,\; t\in [\tau,T],
\end{equation}
where $P$ is the (mild) solution of the following differential Riccati equation
\begin{equation}\label{Riccati}
\addtolength{\jot}{5pt}
\left\{
\begin{alignedat}{2}
&\frac{dP}{dt}-A^*P-PA+PB(K^*K)^{-1}B^*P-C^*C=0\\
&P(0)=0.
\end{alignedat}
\right.
\end{equation}
In other words, $U_T(\cdot,\tau)$ is defined by setting
$$
U_T(t,\tau)x_0:=x(t),\quad \forall x_0\in\calH,\;t\in [\tau,T],
$$
where $x$ is the solution of \eqref{evoopttra}. Here, the solution of $P$ is understood in a mild sense. More precisely, $P\in C_s([0,\infty),\Sigma(\calH))$ is called a mild solution of problem~\eqref{Riccati} if it verifies% the following identity deduced from Duhamel's principle
\begin{equation*}
P(t)x_0=\int_0^t\calT^*_sC^*C\calT_sx_0ds-\int_0^t\calT^*_{t-s}P(s)B(K^*K)^{-1}B^*P(s)\calT_{t-s}x_0ds,
\end{equation*}
for any $x_0\in\calH$ and $t\geq0$.
The next result shows that when the turnpike property holds, the optimal control admits a closed form solution as in \cite[Lemma 3.3.11]{thesis}.
\begin{cor}\label{necess2}
If problem $(GLQ)_T$ satisfies the measure or exponential turnpike property at some steady state $(x_e,u_e)$, then there exists a unique $w\in D(A^*)$ such that
\begin{equation*}
\bvec{A^*}{B^*}w=\bvec{z+C^*Cx_e}{v+K^*Ku_e}.
\end{equation*}
Moreover, the optimal control of problem $(GLQ)_T$ is given in a feedback law form by
\begin{align*}
\begin{split}
u_T^*(t,x_0)-u_e=-(K^*&K)^{-1}B^*P(T-t)(x_T^*(t,x_0)-x_e)\\
&-(K^*K)^{-1}B^*(U_{T-t}(T-t,0))^*w,\;\;\forall t\in [0,T],
\end{split}
\end{align*}
where $P$ is the solution to \eqref{Riccati}.
\end{cor}
\begin{rem}
The vector $w$ defined in the above lemma is called the optimal adjoint steady state of problem $(GLQ)_T$ since it can be seen as the infinite dimensional analogy of the Lagrange multiplier of the optimal steady state problem \eqref{OSSP}.
\end{rem}
\begin{rem}
Following \cite[Part IV, Chapter 1, Proposition 6.2]{Bensou}, the cost functional $J_T$ for problem $(LQ)_T$ can be rewritten as:
\begin{equation}\label{Ptransform}
J_T(x_0,u)=\int_0^T\|K(u(t)+(K^*K)^{-1}B^*P(T-t)x(t))\|^2dt+\<P(T)x_0,x_0\>,
\end{equation}
where $x$ is the solution of \eqref{sy}. In this way, we deduce a decomposition of $J_T$ that explicitly
highlights the dependence on each input $x_0$ and $u$. In particular, the optimal control for problem $(LQ)_T$ is given by: 
$$
u(t)=-(K^*K)^{-1}B^*P(T-t)x(t),\;t\in[0,T].
$$ 
\end{rem}

We can further elaborate on the previous results, and obtain some necessary and sufficient
conditions for the turnpike property in the case of point spectrum. To this aim, let us denote by $\sigma^-(A)$, $\sigma^0(A)$ and $\sigma^+(A)$ the set of all the elements in the spectrum $\sigma(A)$ of $A$ with negative, zero and positive real part, respectively. Following~\cite{Bensou}, we define the point spectrum assumption~\eqref{PS} by
\begin{equation}\label{PS}\tag{$\calP\calS$}
\left\{
\begin{alignedat}{2}
&(a)\;\;\, \text{the set } \sigma^+(A) \text{ consists of a finite set of}\\
&\quad\quad \text{eigenvalues of finite algebraic multiplicity},\\
&(b)\;\;\; \text{there exists } \eps>0, N_A>0 \text{ such that}\\
&\quad\quad \sup_{s\in \sigma^-(A)}\Re s<-\eps,\;\|\calT_t\mathbb P^-_A\|\leq N_Ae^{-\eps t},\;\forall t\geq 0.
\end{alignedat}
\right.
\end{equation}
Here $\mathbb P^-_A$ represents the projector on $\sigma^-(A)$ defined by
\begin{equation*}
\mathbb P^-_A:=\frac{1}{2\pi i}\int_{\gamma^{-}}(sI-A)^{-1}ds
\end{equation*}
where $\gamma^{-}$ is a simple Jordan curve around $\sigma^-(A)$. The projectors $\mathbb P^0_A$ and $\mathbb P^+_A$ are defined analogously. Then $\mathbb P^-_A(\calH)$, $\mathbb P^0_A(\calH)$ and $\mathbb P^+_A(\calH)$ are invariant subspaces for $\calT$.
\begin{rem}\label{PScondition}
By \cite[Part V, Chapter 1, Remark 3.5]{Bensou}, assumptions \eqref{PS} are verified in each of the following cases:
\begin{enumerate}
    \item[$(a)$] $\calH$ is finite dimensional.
    \item[$(b)$] $\calT_t$ is compact for any $t>0$.
\end{enumerate}
\end{rem}

The next theorem provides a complete characterization of the turnpike property in terms of
the stabilizability and detectability of the system for the point spectrum case.

\begin{thm}\label{suffnece1}
If $A$ fulfills assumptions \eqref{PS}, then the following statements are equivalent:
\begin{enumerate}
\item[$(a)$] Problem $(GLQ)_T$ satisfies the exponential turnpike property at some steady state.
\item[$(b)$] Problem $(GLQ)_T$ satisfies the measure turnpike property at some steady state.
\item[$(c)$] The pair $(A,B)$ is exponentially stabilizable and the pair $(A,C)$ is exponentially detectable.
\end{enumerate}
\end{thm}
The following corollary is a direct consequence of the above theorem.
\begin{cor}\label{suffnece2}
If $\calH$, $\calU$ and $\calY$ are all finite dimensional spaces, then the following statements are equivalent:
\begin{enumerate}
\item[$(a)$] Problem $(GLQ)_T$ satisfies the exponential turnpike property at some steady state.
\item[$(b)$] Problem $(GLQ)_T$ satisfies the measure turnpike property at some steady state.
\item[(c)] The pair $(A,B)$ is stabilizable and the pair $(A,C)$ is detectable.
\end{enumerate}
\end{cor}
Our last result shows the exponential turnpike property of the generalized LQ optimal control problem is equivalent to the exponential turnpike property of the LQ optimal control problem.
\begin{thm}\label{suffnece3}
Problem $(GLQ)_T$ satisfies the exponential turnpike property at some steady state $(x_e,u_e)$ if and only if problem $(LQ)_T$
satisfies the exponential turnpike property at $(0,0)$.
\end{thm}

\begin{rem}
A natural question is whether it is possible to extend the above results to the case of unbounded control and observation operators. 

The problem in the unbounded setting quickly becomes far more technical and likely requires a more refined functional analytic framework of the OCP. In general terms, the method of this paper exploits the connection between the differential Riccati equation and the LQ optimal control problem. In particular, we are able to explicitly express the optimal control of~\eqref{OCP} in terms of the solution to~\eqref{Riccati}. However, as noted in~\cite{Weissexample}, in the unbounded setting this becomes a "nasty" problem, since 'various unbounded operators pop up and their domains do
not match’.
%’In the formula linking the optimal feedback operator to the optimal cost operator, as well as in the Riccati equation, the weighting operator of the input has to be replaced by another operator’, ’Despite its simplicity, this is a ”nasty” problem when we want to reconcile it with the existing LQ optimal control theory: various unbounded operators pop up and their domains do
%not match’. 
Similarly, the proper treatment of the unbounded observation case may require the application of different techniques. Even assuming that, using the framework in~\cite[Part IV, Chapter 1, Section 6.2]{Bensou}, we could provide a solution to the differential Riccati equation in the case of admissible unbounded observation, the proof of Theorem~\ref{necess1} heavily relies on the estimates~\eqref{pointcont} and~\eqref{lowerequi2}. However, in the unbounded observation case, $\ell$ would not be defined on $\calH\times\calU$, which makes unclear how to formulate the estimate~\eqref{pointcont}. Besides, since the notion of measure turnpike property does not necessarily control the $L^2$-norm of the optimal pair and $(x_e,u_e)$, it seems difficult to directly deduce~\eqref{lowerequi2}.
\end{rem}

\begin{example}[Parabolic equations]
\label{Ex:ParEqs}
Let $\Omega$ be a bounded domain in $\R^n$ with $\partial\Omega$ of class $C^2$. Let $\calH=\calU=\calY=L^2(\Omega)$, and $B,C:L^2(\Omega)\to L^2(\Omega)$ be linear and bounded operators. Let $K$ be a coercive operator in $\calL(L^2(\Omega))$. We define $A:D(A)\to \calH$ by
\begin{align*}
&D(A)=H^2(\Omega)\cap H_0^1(\Omega),\\
&Ah=(\Delta+cI) h,\quad\forall h\in D(A).
\end{align*}
where $c\in\R$. By \cite[Paty IV, Chapter 1, Section 8.1]{Bensou}, $A$ generates a strongly continuous (in fact, analytic) semigroup on $\calH$. This functional-analytic framework is well suited to describe the distributed
control of the following parabolic equation with Dirichlet boundary conditions:
\begin{equation*}
\left\{
\begin{alignedat}{4}
&\frac{\ptl h}{\ptl t}(x,t)=(\Delta_x+c)h(x,t)+B(u(\cdot,t))(x),\;\;\text{in}\;\;\Omega\times[0,T],\\
&h(\cdot,0)=h_0\in L^2(\Omega),\\
&h(t,x)=0,\;\;\text{on}\;\;\ptl\Omega\times[0,T],
\end{alignedat}
\right.
\end{equation*}
Let $z,v\in L^2(\Omega)$. Consider the optimal control problem: Minimize
\begin{equation*}
J_T(h_0,u)=\int_0^T\int_{\Omega}|(Ch(\cdot,t))(x)-z(x)|^2+|(Ku(\cdot,t))(x)-v(x)|^2dxdt
\end{equation*}
over all $u\in L^2(0,T;L^2(\Omega))$. Obviously this cost functional only differs from the one given in $(GLQ)_T$ by a constant, so our results on turnpike property is valid for this system. 

By \cite[Remark 3.6.4]{TusWei}, $-\Delta$ is a strictly positive operator with compact resolvents, so by spectral theorem, $A$ is diagonalizable with a orthonormal basis $(\psi_k)_{k\in\N}$ of eigenvectors such that the corresponding sequence of eigenvalues $(\lambda_k)_{k\in\N}\subset\R$ is decreasing and satisfies that $\lambda_k\sra-\infty$ as $k\to\infty$. Since (see \cite[Section 2.6]{TusWei}), we can write
\begin{equation*}
\calT_tx=\sum_{k=0}^\infty e^{\lambda_kt}\psi_k,\;\;\forall x\in\calH,t\in[0,\infty),
\end{equation*}
Hence $\calT_t$ is compact for any $t>0$. By Remark \ref{PScondition}, this problem fulfills assumptions \eqref{PS}.

If $\lambda_0<0$, then $\calT$ is exponentially stable, so the pair $(A,B)$ and $(A,C)$ are exponentially stabilizable and exponentially detectable for any $B,C\in \calL(L^2(\Omega))$, respectively. If there exists $m\in\N$ such that $\lambda_m\geq0$ and $\lambda_{m+1}<0$, then by Lemma \ref{Hautus}, the pair $(A,C)$ (resp. $(A,B)$) is exponentially detectable (resp. exponentially stabilizable) if and only if 
\begin{equation}\label{HautType}
\ker (\lambda_iI-A)\cap \ker C=\{0\}\;\;(\text{resp. } \ker (\lambda_iI-A^*)\cap \ker B^*=\{0\}),\;\;\forall i=0,1,\ldots,m.
\end{equation}

As a consequence, the exponential turnpike property is equivalent to the measure turnpike property in this case, and it holds if and only if the above Hautus type conditions hold.
\end{example}
\begin{rem}
There are several common choices of the control operator $B$ and the observation operator $C$ satisfying the requirements of Example~\ref{Ex:ParEqs}. For instance, given an open set $\omega\subset\subset \Omega$, we can consider the following control operators that only act on $\omega$, $B_1:L^2(\omega)\to L^2(\Omega)$, defined by
\begin{equation*}
B_1u(x):=1_{\omega}(x)u(x),\quad\forall x\in \Omega\; ,
\end{equation*}
or $B_2:\R\to L^2(\Omega)$ defined by
\begin{equation*}
B_2u(x):=1_{\omega}(x)u,
\end{equation*}
where $1_{\omega}$ is the characteristic function of $\omega$.

By taking $C_j=B_j^*$, $j = 1,2$, it is not hard to verify that
\begin{equation*}
(C_1h(\cdot))(x)=h(x),\quad\forall h\in L^2(\Omega),\;\forall x\in\omega,
\end{equation*}
and
\begin{equation*}
C_2h(\cdot)=\int_{\omega}h(x)dx,\quad\forall h\in L^2(\Omega)\, .
\end{equation*}

Notice that~\eqref{HautType} is equivalent to whether there exists some $h\in D(A)$ such that, for some $\lambda\geq 0$, we have
\begin{equation}
\left\{
\begin{alignedat}{2}
&(\Delta_x +(c-\lambda)I)h=0,\;\text{a.e.},\\
&h|_{\omega}\equiv0,
\end{alignedat}
\right.
\end{equation}
for the first case and
\begin{equation}\label{Hautussecondcase}
\left\{
\begin{alignedat}{2}
&(\Delta_x +(c-\lambda)I)h=0,\;\text{a.e.},\\
&\int_{\omega}h(x)dx=0,
\end{alignedat}
\right.
\end{equation}
for the second case.

The exponential stabilizability of the pair $(A,B)$ or, equivalently, the exponential detectability of the pair $(A,C)$, for the first case can be deduced from a unique continuation result for second-order elliptic operators~\cite[Theorem 15.2.1]{TusWei}. Namely, if $h(x)\in D(A)$ and $h(x)\equiv0$ for all $x\in \omega$, then $h\equiv 0$ in $\Omega$. For the second case, generally the exponential stabilizability or the exponential detectability may not be guaranteed with arbitrary domains $\Omega$ and $\omega$ even in the $1$-dimensional case unless $c\leq 0$, in which case the semigroup $\calT$ is exponentially stable, and the exponential stabilizability and detectability easily follow. For instance, if we let $\Omega=(0,\pi)$ and $c=5$, then $\lambda=1$ is an eigenvalue of $A=\Delta+cI$, and the corresponding eigenvector is simply $\operatorname{sin}(2x)$. Notice that~\eqref{Hautussecondcase} is now satisfied with $w=(\frac{\pi}{4},\frac{3\pi}{4})$. However, condition~\eqref{Hautussecondcase} provides a sufficient and necessary condition to determine whether the exponential stabilizability or detectability holds. 
\end{rem}

\section{Proofs}\label{Proofs}
\subsection{Proof of Theorem \ref{necess1} (a)}

We first start with some notations.
For each $t\geq 0$, let $\Phi_t:L^2(0,\infty;\calU)\to\calH$ be defined by
\begin{equation*}
\Phi_tu:=\int_0^t\calT_{t-s}Bu(s)ds,
\end{equation*}
and $\Psi_t:\calH\to L^2(0,\infty;\calY)$ be defined by
\begin{equation*}
(\Psi_tx)(s):=\begin{cases}
C\calT_sx,&s\in [0,t],\\
0,&s>t.
\end{cases}
\end{equation*}
Both operators $\Phi_t$ and $\Psi_t$ are bounded, for any $t\ge 0$. In particular, recall that the solution of problem \eqref{sy} is given by
\begin{equation*}
x(t)=\calT_tx_0+\Phi_tu,\;\;\forall t\geq0.
\end{equation*}
The following lemma provides an invariance principle for optimal controls associated with initial conditions in the unobservable subspace $U^\infty$.
\begin{lem}\label{samecontrol}
For any time horizon $T>0$ and $x_0\in U^{\infty}$,
\begin{equation*}
u^*_T(\cdot,x_0)=u^*_T(\cdot,\gamma x_0),\;\;\forall \gamma\in\C.
\end{equation*}
\end{lem}
\begin{proof}
Suppose that $T>0$ and $x_0\in U^{\infty}$. For any $\gamma\in\C$, let $x_\gamma$ be the solution of system \eqref{sy} on $[0,T]$ corresponding to initial condition $\gamma x_0$ and input $u\in L^2(0,T;\calU)$, then
$$
x_\gamma(t)=\calT_t\gamma x_0+\Phi_t(u),\;\;\forall t\in[0,T].
$$
Since $x_0\in U^\infty$, notice that for any $t\in [0,T]$,
\begin{equation*}
\|Cx_\gamma(t)\|^2=\|C\calT_t\gamma x_0\|^2+2\Re\<C\calT_t\gamma x_0,C\Phi_t(u)\>+\|C\Phi_t(u)\|^2=\|C\Phi_t(u)\|^2,
\end{equation*}
and
\begin{equation*}
2\Re\<z,x_\gamma(t)\>=2\Re\<z,\calT_t\gamma x_0\>+2\Re\<z,\Phi_t(u)\>.
\end{equation*}
So, we have
\begin{align*}
\begin{split}
J_T(\gamma x_0,u)
=&\int_0^T\|Cx_\gamma(t)\|^2 + \|Ku(t)\|^2 + 2\Re\< z,x(t)\> + 2\Re\< v,u(t)\> dt\\
=&\int_0^T \|C\Phi_t(u)\|^2+2\Re\< z,\calT_t\gamma x_0+\Phi_t(u)\> + \|Ku(t)\|^2 + 2\Re\<v,u(t)\>dt.
\end{split}
\end{align*}
We thus realize that the terms in $J_T(\gamma x_0,u)$ containing the input $u$ do not depend on the choice of $\gamma\in\C$. Hence, we conclude that $u^*_T(\cdot,\gamma x_0)=u^*_T(\cdot,x_0)$ for any $\gamma\in \C$.
\end{proof}

\begin{proof}[Proof of Theorem \ref{necess1} (a)]

Without loss of generality, assume that the measure turnpike property is satisfied at $(x_e,u_e)=(0,0)$. Otherwise, by taking $\wt x:=x-x_e$ and $\wt u=u-u_e$, it is not hard to verify that for any $x_0\in\calH$, if $(x^*(\cdot),u^*(\cdot))$ minimizes the cost functional \eqref{OCP}, then $(\wt x^*(\cdot),\wt u^*(\cdot)):=(x^*(\cdot)-x_e,u^*(\cdot)-u_e)$ minimizes the cost functional
\begin{equation}\label{OCPchangevariable}
J_T(\wt x_0,\wt u):=\int_0^T\wt \ell(\wt x(t),\wt u(t))dt,
\end{equation}
over any pair $(\wt x(\cdot),\wt u(\cdot))\in L^2(0,T;\calH)\times L^2(0,T;\calU)$ constrained to
\begin{equation}\label{sychangevariable}
\dot {\wt x}(t)=\dot x(t)-x_e=Ax(t)+Bu(t)=Ax(t)+Bu(t)-Ax_e-Bu_e=A\wt x(t)+B\wt u(t),
\end{equation}
where $\wt x_0:=x_0-x_e$ and
\begin{equation}\label{costchangevariable}
\wt{\ell}(\wt x,\wt u) := \|C\wt x\|^2 + \|K\wt u\|^2 + 2\Re\langle z+C^*Cx_e,\wt x\rangle + 2\Re\langle v+K^*Ku_e,\wt u\rangle.
\end{equation}
By Definition~\ref{turnpikedef}, the modified OCP~\eqref{OCPchangevariable}-\eqref{costchangevariable} satisfies the measure turnpike property at $(0,0)$, whereas the pair $(A,C)$ is unchanged; thus, Theorem~\ref{necess1} (a) also remains unchanged.

Now, let $\calN$ be the closed unit ball centered at $0$ in $\calH$. Observe that there exists a sufficiently large constant $M>0$ such that, if $x_0\in \calH$ and $t_0\geq 2M_{\calN,\,\eps}+1$ where $M_{\calN,\,\eps}$ is defined as in Definition \ref{turnpikedef}, then
\begin{equation*}
\|\calT_{t_0}x_0\|\geq M \implies\|\calT_{t}x_0\|>2\eps,\;\;\forall t\in [t_0-2M_{\calN,\,\eps}-1,t_0].
\end{equation*}
{In fact, notice that there exist some positive constants $k$ and $M_k$ (see, e.g., \cite[Proposition 2.1.2]{TusWei}) such that
\begin{multline*}
M\leq \|\calT_{t_0}x_0\|\leq \|\calT_{t_0-t}\|\|\calT_{t}x_0\| \leq M_k\, e^{k(t_0-t)}\|\calT_{t}x_0\|\\
\le M_k\, e^{k(2M_{\calN,\,\eps}+1)}\|\calT_{t}x_0\|\,,\quad\forall t\in [t_0-2M_{\calN,\,\eps}-1,t_0].
\end{multline*}
Thus we may conclude that, if $M>2\eps M_k\, e^{k(2M_{\calN,\,\eps}+1)}$, then $\|\calT_{t}x_0\|>2\eps$.}

We now argue by contradiction, and assume that $\calT$ is not exponentially stable on $U^\infty$. Hence, the restriction of $\calT$ on \(U^{\infty}\) does not converge to $0$ in operator norm as \(t\to \infty\)~\cite[Proposition~V1.2]{Oneparametersemigroup}. In other words, there exists some $\eps>0$ such that
\begin{equation*}
\limsup_{t\geq 0}\|\calT_t|_{U^{\infty}}\|>\eps.
\end{equation*}
As a consequence, we are able to find some $t_r>2M_{\calN,\,\frac{\eps}{4M}}+2M_{\calN,\,\eps}+2$ and $x_0\in U^\infty$ such that $\|x_0\|\leq 1$ and $\|\calT_{t_r}x_0\|> \eps$.

By measure turnpike property, there exists a point $t_l\in[0,2M_{\calN,\,\frac{\eps}{4M}}+1]$ such that
\begin{equation}\label{bothinneibor}
\Big\|x^*_{t_r}(t_l,x_0)\Big\|\leq \frac{\eps}{4M}\;\;\text{and}\;\;\left\|x^*_{t_r}\left(t_l,\frac{x_0}{2}\right)\right\|\leq\frac{\eps}{4M}.
\end{equation}
In fact, since $x_0,\frac{x_0}{2}\in\calN$, the measure of the sets of $t\in[0,t_r]$ such that $\|x^*_{t_r}(t,x_0)\|>\frac{\eps}{4M}$ and $\|x^*_{t_r}(t,\frac{x_0}{2})\|>\frac{\eps}{4M}$ are both less than or equal to $M_{\calN,\,\frac{\eps}{4M}}$, but $[0,t_r]$ has a measure larger than $2M_{\calN,\,\frac{\eps}{4M}}$, so there must be some point in $[0,t_r]$ that satisfies \eqref{bothinneibor}. 

Since $x_0\in U^\infty$, from Lemma \ref{samecontrol}, we obtain that
\begin{align*}
\frac{\|\calT_{t_l}x_0\|}{2}=\left\|x^*_{t_r}(t_{l},x_0)-x^*_{t_r}\left(t_{l},\frac{x_0}{2}\right)\right\|\leq \Big\|x^*_{t_r}(t_l,x_0)\Big\|+\left\|x^*_{t_r}\left(t_l,\frac{x_0}{2}\right)\right\|\leq \frac{\eps}{2M}.
\end{align*}
This implies that $\|\frac{M}{\eps}\calT_{t_l}x_0\|\leq 1$. Thus $\frac{M}{\eps}\calT_{t_l}x_0\in \calN$. Besides, notice that
\begin{equation*}
\left\|\calT_{t_r-t_l}\frac{M}{\eps}\calT_{t_l}x_0\right\|=\left\|\frac{M}{\eps}\calT_{t_r}x_0\right\|>\frac{M}{\eps}\eps\geq M\;\;\text{and}\;\; t_r-t_l>2M_{\calN,\,\eps}+1,
\end{equation*}
so
\begin{equation*}
\left\|\calT_{t}\frac{M}{\eps}\calT_{t_l}x_0\right\|>2\eps,\;\;\forall t\in [t_r-t_l-2M_{\calN,\,\frac{\eps}{2}}-1,t_r-t_l].
\end{equation*}
Applying again Lemma \ref{samecontrol}, we obtain that
\begin{equation*}
\left\|x^*_{t_r}\left(t,\frac{M}{\eps}\calT_{t_l}x_0\right)-x^*_{t_r}\left(t,\frac{1}{2}\frac{M}{\eps}\calT_{t_l}x_0\right)\right\|=\left\|\frac{1}{2}\calT_{t}\frac{M}{\eps}\calT_{t_l}x_0\right\|>\eps
\end{equation*}
holds for any $t\in[t_r-t_l-2M_{\calN,\,\frac{\eps}{2}}-1,t_r-t_l]$. However, this implies a contradiction: indeed, the length of this interval is greater than or equal to $2M_{\calN,\,\frac{\eps}{2}}+1$, and similar arguments as above show that there must exist a $t\in [t_r-r_l-2M_{\calN,\,\frac{\eps}{2}}-1,t_r-t_l]$ such that
\begin{equation*}
\left\|x^*_{t_r}\left(t,\frac{M}{\eps}\calT_{t_l}x_0\right)\right\|\leq \frac{\eps}{2}\;\;\text{and}\;\;\left\|x^*_{t_r}\left(t,\frac{1}{2}\frac{M}{\eps}\calT_{t_l}x_0\right)\right\|\leq \frac{\eps}{2},
\end{equation*}
which further implies
\begin{equation*}
\left\|x^*_{t_r}\left(t,\frac{M}{\eps}\calT_{t_l}x_0\right)-x^*_{t_r}\left(t,\frac{1}{2}\frac{M}{\eps}\calT_{t_l}x_0\right)\right\|\leq \frac{\eps}{2}+\frac{\eps}{2}\leq \eps.
\end{equation*}
\end{proof}

\subsection{Proof of Theorem \ref{necess1} (b)}
We first prove two preliminary lemmas.

\begin{lem}[Characterization of the optimal steady state]\label{uniqueoss}
$(x_e,u_e)\in D(A)\times\calU$ is an optimal steady state of the optimal steady state problem \eqref{OSSP} if and only if
\begin{equation*}
\bvec{z+C^*Cx_e}{v+K^*Ku_e}\in \ol{\ran}\bvec{A^*}{B^*}.
\end{equation*}
Moreover, the optimal steady state is unique if and only if
\begin{equation*}
\ker A\cap \ker C=\{0\}.
\end{equation*}
\end{lem}
\begin{proof}
Let $V=\ker[A\;B]$. Observe that $\ell$ restricted to $V$ can be written as
$$
\ell\left(x,u\right)=\left\<\mathbb{P}_V \bmat{C^{*}C}{0}{0}{K^{*}K}\bigg|_{V}\bvec{x}{u}+2\mathbb{P}_{V}\bvec{z}{v},\bvec{x}{u}\right\>,
$$
where $\mathbb{P}_{V}$ is the projection vector onto the closed subspace $V$.

Let the non-negative operator $\calP:V\to V$ be defined by 
$$
\calP:=\mathbb{P}_{V} \bmat{C^{*}C}{0}{0}{K^{*}K}\bigg|_{V}.
$$
Notice that $\calP\in L(V)$ and $\mathbb{P}_{V}\bvec{z}{v}\in V$. So, by \cite[Lemma 4]{LFM}, $(x_e,u_e)\in V$ is a minimizer of problem \eqref{OSSP} if and only if
\begin{equation}\label{chaross}
\mathbb{P}_{V}\bmat{C^{*}C}{0}{0}{K^{*}K}\bigg|_{V}\bvec{x_e}{u_e}=-\mathbb{P}_{V}\bvec{z}{v}.
\end{equation}
Equivalently, $(x_e,u_e)\in V$ satisfies
$$
\bvec{z+C^*Cx_e}{v+K^*Ku_e}=\bmat{C^{*}C}{0}{0}{K^{*}K}\bigg|_{V}\bvec{x_e}{u_e}+\bvec{z}{v}\in V^{\perp}=\ol{\ran}\bvec{A^*}{B^*}.
$$
This proves our first claim.

On the other hand, simple considerations show that the uniqueness of $(x_e,u_e)\in V$ verifying equation \eqref{chaross} is equivalent to $\ker \calP=\{(0,0)\} \subset D(A)\times\calU$. Now let us show that $\ker \calP=\{(0,0)\}$ is further equivalent to $\ker A\cap \ker C=\{0\}$.

Suppose that $\ker A\cap \ker C=\{0\}$. If $(x_0,u_0)\in \ker \calP$, then we have that
$$
0 = \left\<\calP\bvec{x_0}{u_0},\bvec{x_0}{u_0}\right\>=\left\<\bvec{C^{*}Cx_0}{K^{*}Ku_0},\bvec{x_0}{u_0}\right\>=\|Cx_0\|^2+\|Ku_0\|^2\; .
$$
Since $K$ is coercive, we deduce that $u_0=0$ and $Cx_0=0$. Besides, since $(x_0,u_0)\in V$, we have that $Ax_0 + Bu_0 = 0$. However, since $u_0 = 0$, this implies that $Ax_0=-Bu_0=0$. Hence, $Ax_0 = Cx_0 = 0$. This implies that $x_0\in \ker A\cap \ker C$, which by assumption is trivial. Thus $x_0=0$ and $\ker \calP=\{(0,0)\}$. 

Conversely, if there exists some $x_0\in \ker A\cap \ker C$ and $x_0\neq 0$, we can easily verify that $(x_0,0)\in \ker\calP\subset V$. So, $\ker \calP\neq\{(0,0)\}$.

\end{proof}

The following technical lemma provides a very conservative estimate for the lower bound of the cost functional.
\begin{lem}\label{lowbound}
Suppose $\calN$ is a bounded subset of $\calH$. Then for any $T\in (0,\infty)$ and $x_0\in\calN$, 
$$
\min_{u\in L^2(0,T;\mathcal{U})}J(x_0,u)\geq -M(e^{2kT}-1)
$$
for some constants $M,k>0$.
\end{lem}
\begin{proof}
Since $K$ is coercive, $K^*K$ is bounded below and thus invertible. Let $u_0:=-(K^*K)^{-1}v$, then there exists $M_K>0$ such that
\begin{equation}\label{ubelow}
\|Ku\|^2 + 2\Re\<v,u\>=\|K(u-u_0)\|^2-\|Ku_0\|^2\geq M_K\|u-u_0\|^2-\|Ku_0\|^2,\;\;\forall u\in\calU.
\end{equation}
Besides, by \cite[Proposition 2.1.2 and Proposition 4.3.3]{TusWei}, there exist positive constants $k$, $M_1$ and $M_2$ such that
\begin{equation*}
\|\Phi_t\|\leq M_1e^{kt},\;\; \text{and}\;\; \|\calT_t\|\leq M_2e^{kt},\;\;\forall t\geq0
\end{equation*}

Assume that $T>0$, and let $x$ be the solution of \eqref{sy} corresponding to the initial condition $x_0\in \calN$ and input $u\in L^2(0,T;\calU)$. Then, for any $t\in [0,T]$,
\begin{align}\label{xbelow}
\begin{split}
2\Re\< z,x(t)\>&\geq -2\|z\|\left(\|\calT_tx_0\|+\|\Phi_t(u-u_0)\|+\|\Phi_tu_0\|\right)\\
&\geq -2\|z\|(M_2e^{kt}\|x_0\|+M_1e^{kt}\|u-u_0\|_{L^2}+\|\Phi_tu_0\|).
\end{split}
\end{align}
Concerning the term with the norm of $\Phi_tu_0$, we have that
\begin{align}\label{phiu_0}
\begin{split}
\|\Phi_tu_0\|&=\left\|\int_0^t\calT_{t-s}Bu_0ds\right\|\leq \int_0^tM_2e^{k(t-s)}\|Bu_0\|ds\leq \frac{e^{kt}}{k}M_2\|Bu_0\|.
\end{split}
\end{align}

Now, \eqref{ubelow}, \eqref{xbelow} and \eqref{phiu_0} together imply that there exist some positive constants $M_3$ and $M_4$ such that
\begin{equation*}
\ell(x(t),u(t))\geq -M_3e^{kt}-M_4e^{kt}\|u-u_0\|_{L^2}+M_K\|u(t)-u_0\|^2-\|Ku_0\|^2,\;\;\forall t\in[0,T].
\end{equation*}
Integrating $\ell$ over $[0,T]$, a straightforward calculation gives
\begin{align*}
\begin{split}
J_T(x_0,u)\geq -\frac{M_4}{k}(e^{kT}-1)\|u-u_0\|_{L^2}-\frac{M_3}{k}(e^{kT}-1)+M_K\|u-u_0\|^2_{L^2}-\|Ku_0\|^2T.
\end{split}
\end{align*}
Besides, notice that, as a quadratic function with respect to $\|u-u_0\|_{L^2}$,
\begin{equation*}
M_K\|u-u_0\|^2_{L^2}-\frac{M_4}{k}(e^{kT}-1)\|u-u_0\|_{L^2}\geq-\frac{M_4^2(e^{kT}-1)^2}{4M_0k^2}.
\end{equation*}

Simple considerations about the above two inequalities show that
\begin{equation*}
J_T(x_0,u)\geq -M_5(e^{2kT}-1)-M_6(e^{kT}-1)-M_7T
\end{equation*}
holds for some positive constants $M_5$, $M_6$ and $M_7$. Since the dominant term of the above estimate is $e^{2kT}-1$, there exists some $M>0$ (independent of $T$) such that
$$
J_T(x_0,u)\geq -M(e^{2kT}-1).
$$
The lemma then follows from the fact that $u\in L^2(0,T;\calU)$ can be arbitrarily chosen.

\end{proof}

Now we are in the position to prove Theorem \ref{necess1} (b).
\begin{proof}[Proof of Theorem \ref{necess1} (b)]
We first show that if the measure turnpike property holds at some steady state $(x_e,u_e)$, then $(x_e,u_e)$ is necessarily an optimal steady state. 

We prove this by contradiction. Suppose there exists some steady state $(\wt{x}_e,\wt{u}_e)$ such that $\ell(\wt x_e,\wt u_e)<\ell(x_e,u_e)$. Define $(dx,du):=(\wt x_e,\wt u_e)-(x_e,u_e)$. Fix some bounded neighborhood $\calN$ of $x_e$, and some sufficiently small $\la\in (0,1]$ such that $x_e+\la dx\in \calN$. Then, by the convexity of $\ell$,
$$
\ell(x_e+\la dx,u_e+\la du)<\ell(x_e,u_e).
$$ 
Now let us denote the new equilibrium point $(x_e+\la dx,u_e+\la du)$ by $(\wt x_e,\wt u_e)$ again. 

Since $\ell$ is continuous, there exist some sufficiently small $\eps$ and $\delta>0$ so that, for any $(x_0,u_0)\in \calH\times\calU$, 
\begin{equation}\label{pointcont}
\|x_0-x_e\|+\|u_0-u_e\|\leq\delta \implies \ell(\wt x_e,\wt u_e)+\eps<\ell(x_0,u_0).
\end{equation}
Now fix some $T>0$. Then we know that the set
$$
A:=\{t\in[0,T]\;\big|\; \|u^*_T(t,\wt x_e)-u_e\|+\|x^*_T(t,\wt x_e)-x_e\|>\delta\}
$$
is open, and its Lebesgue measure is smaller than $M_{\calN,\delta}$, where $M_{\calN,\delta}$ is defined as in Definition~\ref{turnpikedef}.  

Without loss of generality, assume that $A=\bigcup_{j=1}^{\infty}(t_{j,l},t_{j,r})$ where $((t_{j,l},t_{j,r}))_{j\in \N}$ are disjoint open intervals (the number of such intervals may be finite, and there may be at most two intervals that contain an endpoint $0$ or $T$, but the proof is basically the same). 

Observe that at each left endpoint $t_{j,l}$, $j\in\N$, of these intervals, we must have 
$$
\|x^*_T(t_{j,l},\wt x_e)-\wt x_e\|\leq \|x^*_T(t_{j,l},\wt x_e)-\wt x_e\|+\|u^*_T(t_{j,l},\wt x_e)-\wt u_e\| =\delta.
$$
From Lemma~\ref{lowbound}, there exist some $M,k>0$ such that
\begin{equation*}
-M(e^{k(t_{j,r}-t_{j,l})}-1)\leq \int_{t_{j,l}}^{t_{j,r}}\ell(x^*_T(t,\wt x_e),u^*_T(t,\wt x_e))dt.
\end{equation*}
Owing to the countable additivity of the Lebesgue integral, we have
\begin{equation*}
-M\sum_{j=1}^{\infty}(e^{k(t_{j,r}-t_{j,l})}-1) \leq \int_{A}\ell(x^*_T(t,\wt x_e),u^*_T(t,\wt x_e))dt.
\end{equation*}
Now simple calculation shows that, for any $t_1$ and $t_2>0$,
\begin{equation*}
-M(e^{k(t_1+t_2)}-1)\leq -M(e^{kt_1}-1)-M(e^{kt_2}-1).
\end{equation*}
Combining the two inequalities, we obtain that
\begin{equation}\label{lowerequi1}
-M(e^{kM_{\calN,\delta}}-1)\leq -M(e^{k\mu(A)}-1)\leq \int_{A}\ell(x^*_T(t,\wt x_e),u^*_T(t,\wt x_e))dt.
\end{equation}
Meanwhile, for any $t\in [0,T]\setminus A$, since $\|x^*_T(t,\wt x_e)-\wt x_e\|+\|u^*_T(t,\wt x_e)-\wt u_e\|\leq\delta$, we have
$$
 \ell(\wt x_e,\wt u_e)+\eps<\ell(x^*_T(t,\wt x_e),u^*_T(t,\wt x_e)),
$$
so
\begin{equation}\label{lowerequi2}
(T-M_{\calN,\delta})(\ell(\wt x_e,\wt u_e)+\eps)< \int_{[0,T]\setminus A}\ell(x^*_T(t,\wt x_e),u^*_T(t,\wt x_e))dt.
\end{equation}
Finally, since $(\wt x_e,\wt u_e)$ is a steady state,
\begin{equation}\label{upperequi}
\int_0^T\ell(x^*_T(t,\wt x_e),u^*_T(t,\wt x_e))dt\leq J_T(\wt x_e,\wt u_e)=T \ell(\wt x_e,\wt u_e).
\end{equation}
Combining \eqref{lowerequi1}, \eqref{lowerequi2} and \eqref{upperequi}, simple calculations show that
\begin{equation*}
T\eps<M(e^{kM_{\calN,\delta}}-1)+M_{\calN,\delta}(\ell(\wt x_e,\wt u_e)+\eps).
\end{equation*}
This leads to a contradiction when $T$ is taken sufficiently large, since the right-hand side does not depend on $T$. We have thus proved that $(x_e,u_e)$ is necessarily an optimal steady state. 

We shall now prove the uniqueness.
Recall from Lemma \ref{uniqueoss} that the uniqueness of the optimal steady state is characterized by 
\begin{equation*}
\ker A\cap \ker C=\{0\}.
\end{equation*}
In fact, if $x\in \ker A\cap \ker C$, {then $(\calT_tx)'=\calT_tAx=0$, and thus $\calT_tx=\calT_0x=x$, $\forall t\geq0$.} So, $C\calT_tx=C\calT_0x=0$ for any $t\geq0$. By definition, $x$ belongs to the unobservable subspace $U^\infty$ of $\calT$. From Theorem~\ref{necess1} (a), we know that $\calT$ is exponentially stable on $U^\infty$. We thus deduce that $\|\calT_tx\|=\|x\|\to0$ as $t\to\infty$. Thus $x=0$.
\end{proof}

\subsection{Proof of Theorem \ref{necess1} (c)}
Define the cost functional $J_{\infty}: \calH\times L^2(0,\infty;\calU)\to \ol \R$ by
\begin{equation*}
J_{\infty}(x_0,u):=\int_0^{\infty}\|Cx(t)\|^2+\|Ku(t)\|^2dt,
\end{equation*}
where $x$ is the solution of \eqref{sy}. Now, we recall from \cite[Part V, Chapter 1, Remark 3.2]{Bensou} an equivalent formulation of exponential stabilizability.
\begin{prop}\label{stabilizabilitychar}
The pair $(A,B)$ is exponentially stabilizable if and only if
\begin{equation}\label{Istabilizability}
\inf_{u\in L^2(0,\infty;\calH)}J_{\infty}(x_0,u)<\infty.
\end{equation}
holds with $C=I$ and any $x_0\in\calH$.
\end{prop}
The following lemma provides a conservative estimate on the upper bound of the $L^2$-norm of the optimal control when the initial point lies in a bounded set. 
%{\ZL Recall that we abbreviate the $L^2$-norm considered on interval $(0,T)$ as $\|\cdot\|_{L^2}$.}
\begin{lem}\label{normonebound}
Assume that $\calN$ is a bounded set in $\calH$. Then for any $T>0$, there exists some $M_u>0$ (dependent on $T$) such that
\begin{equation*}
\sup_{x_0\in\calN}\|u^*_T(\cdot,x_0)\|_{L^2}\leq M_u.
\end{equation*}
\end{lem}
\begin{proof}
Fix some $T>0$. We first prove that there exists some $M_{0}>0$ such that
\begin{equation}\label{uniboundu}
\sup_{x_0\in\calN}\int_0^T\ell(x^*_T(t,x_0),u^*_T(t,x_0))dt\leq M_0.
\end{equation}
To see this, consider the case of input $u\equiv0$, then we obtain
\begin{equation*}
\int_0^T\ell(x^*_T(t,x_0),u^*_T(t,x_0))dt \leq\int_0^T\|C\calT_tx_0\|^2+2\Re\<z,\calT_tx_0\>dt,\;\;\forall x_0\in\calN.
\end{equation*}
Since $\calT_{t}x_0$ is bounded in norm on $[0,T]$ for any $x_0\in\calN$, simple considerations show that the positive constant $M_{0}$ exists.

Now notice that, for any $x_0\in\calN$ and $t\in[0,T]$,
\begin{equation*}
2\Re\<z, x^*_T(t,x_0)\>\geq-2\|z\|\|\calT_tx_0+\Phi_tu^*_T(\cdot,x_0)\|\geq -M_1-M_2\|u^*_T(\cdot,x_0)\|_{L^2}
\end{equation*}
holds for some appropriate $M_1$, $M_2>0$. Therefore,
\begin{align*}
\int_0^T\ell(x^*_T(t,x_0),u^*_T(t,x_0))dt&\geq \int_0^T-M_1-M_2\|u^*_T(\cdot,x_0)\|_{L^2}\\
&\qquad\qquad+M_K\|u^*_T(t,x_0)\|^2-2\|v\|\|u^*_T(t,x_0)\|dt\\
&\geq -M_1T-M_2T\|u^*_T(\cdot,x_0)\|_{L^2}+M_K\|u^*_T(\cdot,x_0)\|^2_{L^2}\\
&\qquad\qquad-2\|v\|\|u^*_T(\cdot,x_0)\|_{L^1}
\end{align*}
holds with some $M_K,M_3>0$, where $\|\cdot\|_{L^1}$ denotes the $L^1$-norm on $L^1(0,T;\calU)$.
\\
By Hölder's inequality, the $L^1$-norm is dominated by the $L^2$-norm when $T$ is fixed, so there exists $M_4>0$ such that
\begin{equation*}
\int_0^T\ell(x^*_T(t,x_0),u^*_T(t,x_0))dt\geq -M_1T-M_4\|u^*_T(\cdot,x_0)\|_{L^2}+M_K\|u^*_T(\cdot,x_0)\|^2_{L^2}.
\end{equation*}

Now simple considerations about equation \eqref{uniboundu} and the above two inequalities show that the $L^2$-norm of $u^*_T(\cdot,x_0)$ is uniformly bounded above for all $x_0\in\calN$, i.e., there exists some positive constant $M_{u}$ (dependent on $T$) such that
\begin{equation*}
\sup_{x_0\in\calN}\|u^*_T(\cdot,x_0)\|_{L^2}\leq M_u.
\end{equation*}
\end{proof}

\begin{proof}[Proof of Theorem \ref{necess1} (c)]
Assume the measure turnpike property is satisfied at $(x_e,u_e)$. Without loss of generality, let $\calN$ be the closed unit ball in $\calH$ centered at $x_e$. 

We claim that, for any $x_0\in\calH$, there exists some $u\in L^2(0,\infty;\calU)$ such that the corresponding trajectory $x$ with input $u$ and initial condition $x_0$ is $L^2$-integrable (on $(0,\infty)$). The case for $x_0=0$ is trivial, so let us assume that $x_0\neq0$.

Fix some $T>M_{\calN,\frac{1}{2}}+2${, where $M_{\calN,\frac{1}{2}}$ is defined as in Definition \ref{turnpikedef}.} Let $t_0:=0$ and $\wt x_0:=\frac{x_0}{\|x_0\|}+x_e\in\calN$. Notice that
\begin{equation*}
\mu\{t\in[0,T]\,|\, \|x^*_T(t,\wt x_0)-x_e\|+\|u^*_T(t,\wt x_0)-u_e\|>1/2\}\leq M_{\calN,\frac{1}{2}}.
\end{equation*}
Since $T>M_{\calN,\frac{1}{2}}+2$, there exists a $t_1>1$ satisfying
\begin{equation*}
\|x^*_T(t_1,\wt x_0)-x_e\|\leq \frac{1}{2}.
\end{equation*}
We define $u(t)$ on $[t_0,t_1)$ by setting
\begin{equation*}
u(t):=\|x_0\|\left(u^*_T(t,\wt x_0)-u_e\right),\;\forall t\in[t_0,t_1).
\end{equation*}

We claim that the solution of system \eqref{sy} on $[t_0,t_1]$ corresponding to initial condition $x_0$ and input $u$, denoted by $x$, is given by
\begin{equation*}
x(t)=\|x_0\|\left(x^*_T(t,\wt x_0)-x_e\right),\;\;\forall t\in[t_0,t_1].
\end{equation*}
In fact, we can verify this by noticing that $x(0)=x_0$ and
\begin{align*}
\dot{x}(t)&=\|x_0\|(Ax^*_T(t,\wt x_0)+Bu^*_T(t,\wt x_0))\\
&=Ax(t)+A\|x_0\|x_e+Bu(t)+B\|x_0\|u_e\\
&=Ax(t)+Bu(t).
\end{align*}
So, we have
\begin{equation*}
\|x(t_1)\|\leq \frac{1}{2}\|x_0\|.
\end{equation*}

On the other hand, notice that
\begin{align*}
\int_{t_0}^{t_1}\|x(t)\|^2+\|u(t)\|^2dt&\leq \|x_0\|^2(T\|x^*_T(\cdot,\wt x_0)\|^2_{L^\infty}+2T\|x_e\|\|x^*_T(\cdot,\wt x_0)\|_{L^\infty}+T\|x_e\|^2\\
&\qquad\quad+\|u^*_T(\cdot,\wt x_0)\|^2_{L^2}+2\|u_e\|\|u^*_T(\cdot,\wt x_0)\|_{L^1}+T\|u_e\|^2),
\end{align*}
where the $L^1$, $L^2$ and $L^\infty$-norms here are considered on the domain $[t_0,t_1]$. Recall from Lemma~\ref{normonebound} that there exists some $M_1>0$ such that
\begin{equation*}
\|u^*_T(\cdot,\wt x_0)\|_{L^2}<M_1,\;\;\forall \wt x_0\in\calN.
\end{equation*}
By Hölder's inequality, the $L^2$-norm dominates the $L^1$-norm, so there exists some $M_2>0$ such that
\begin{equation*}
\|u^*_T(\cdot,\wt x_0)\|_{L^1}\leq M_2\;\;\forall \wt x_0\in\calN.
\end{equation*}
Finally, since {$\Phi_t\in\calL(L^2(0,T;\calU),\calH)$,}
\begin{equation*}
\|x^*_T(t,\wt x_0)\|\leq \|\calT_t\wt x_0\|+\|\Phi_t\|\|u^*_T(\cdot,\wt x_0)\|_{L^2},
\end{equation*}
hence $x^*_T(t,\wt x_0)$ is uniformly bounded in norm for any $\wt x_0\in\calN$ and $t\in[t_0,t_1]$. That is, there exists some positive constant $M_3$ such that
\begin{equation*}
\|x^*_T(t,\wt x_0)\|\leq M_3,\;\;\forall \wt x_0\in\calN,t\in[t_0,t_1].
\end{equation*}
Combining the previous estimates, we deduce that there exists $M>0$ such that
\begin{equation*}
\int_{t_0}^{t_1}\|x(t)\|^2+\|u(t)\|^2dt<\|x_0\|^2M,\;\;\forall\wt x_0\in \calN.
\end{equation*}
It should be noted that we can take the above $M$ in a way that it only depends on the selection of $T$, but not on the selection of the length of the interval $t_1-t_0$.
Now we can repeat the above argument with $t_1$ and $x(t_1) = x_1$ in place of $t_0$ and $x(t_0)=x_0$, respectively. By induction, we can find a sequence $(t_i)_{i\in\N}$ satisfying, and let $u$ be defined inductively on each interval $[t_i,t_{i+1})$ by setting
\begin{equation*}
u|_{[t_i,t_{i+1})}(t):=\|x(t_i)\|\left(u^*_T(t,\wt x_i)-u_e\right)
\end{equation*}
where $\wt x_i:=\frac{x(t_i)}{\|x(t_i)\|}+x_e$, $i=0,1,2,\ldots\,$. In particular, our construction of $t_i$ gives
\begin{equation*}
\|x(t_i)\|\leq \left(\frac{1}{2}\right)^{i}\|x_0\|,\;\;\forall i\in\N.
\end{equation*}
Finally, since for each $i\in\N$, $\wt x_i$ belongs to $\calN$, we have that
\begin{equation*}
\int_{t_i}^{t_{i+1}}\|x(t)\|^2+\|u(t)\|^2dt<\|x(t_i)\|^2M\leq \left(\frac{1}{4}\right)^{i}\|x_0\|^2M.
\end{equation*}
So,
\begin{equation*}
\int_{0}^{\infty}\|x(t)\|^2+\|u(t)\|^2dt\leq \frac{4}{3}\|x_0\|^2M.
\end{equation*}
Theorem \ref{necess1} (c) now follows easily from Proposition \ref{stabilizabilitychar}.
\end{proof}

\subsection{Proof of Corollary \ref{necess2}}

We recall the following result, proved in \cite[Lemma 3.3.11]{thesis}.
\begin{lem}\label{test}
If the pair $(A,B)$ is exponentially stabilizable, then $\ran[A\;B]=\calH$.
\end{lem}

\begin{proof}[Proof of Corollary \ref{necess2}]
Assume that problem $(GLQ)_T$ satisfies the measure turnpike property at some steady state $(x_e,u_e)$. Then, by Lemma \ref{uniqueoss} and Theorem \ref{necess1} (b), we have that
\begin{equation*}
\bvec{z+C^*Cx_e}{v+K^*Ku_e}\in \ol{\ran}\bvec{A^*}{B^*}.
\end{equation*}
Besides, by Theorem \ref{necess1} (c), Lemma \ref{test} and the closed range theorem, $\ran\bvec{A^*}{B^*}$ is closed and
\begin{equation*}
\ker\bvec{A^*}{B^*}=({\ran}[A\;B])^{\perp}=\calH^{\perp}=\{0\}.
\end{equation*}
So, there exists a unique $w\in D(A^*)$ satisfying
\begin{equation*}
\bvec{A^*}{B^*}w=\bvec{z+C^*Cx_e}{v+K^*Ku_e}.
\end{equation*}

Following line by line the
proof of Lemma \ref{explicitoc}, we can show that the optimal control of problem $(GLQ)_T$ is given in a feedback law form by
\begin{align*}
\begin{split}
u_T^*(t,x_0)-u_e=-(K^*&K)^{-1}B^*P(T-t)(x_T^*(t,x_0)-x_e)\\
&-(K^*K)^{-1}B^*(U_{T-t}(T-t,0))^*w,\;\;\forall t\in [0,T],
\end{split}
\end{align*}
where $P$ is the solution to~\eqref{Riccati}.
\end{proof}
\begin{rem}
It should be noted that the explicit formula of the optimal control provided in Corollary~\ref{necess2} was first obtained in~\cite{Zua} in a slightly different framework. In our case, the proof of this formula crucially relies on Lemma~\ref{test}.
\end{rem}

\subsection{Proof of Theorem \ref{suffnece1}}
The following lemma provides a Hautus type condition for detectability~\cite[Part V, Chapter 1, Proposition 3.3]{Bensou}, which will be useful in the following proof.
\begin{lem}\label{Hautus}
If $A$ fulfills assumptions \eqref{PS}, then the following conditions are equivalent:
\begin{enumerate}
\item[$(a)$] The pair $(A,C)$ is exponentially detectable.
\item[$(b)$] $\ker(sI-A)\cap\ker C=\{0\}\ $ for all $s\in \sigma^0(A)\cup\sigma^+(A)$.
\end{enumerate}
\end{lem}
%Next, let us provide a proof of Theorem \ref{suffnece1}.

\begin{proof}[Proof of Theorem \ref{suffnece1}]
$(a)\Sra(b)$: This is trivial. $(c)\Sra(a)$: This result has been proved in several different papers, of which we refer to \cite[Theorem 3.3.1]{thesis} and \cite[Theorem 3]{TreZua}.

So, we only need to verify that $(b)\Sra(c)$. Suppose that problem $(GLQ)_T$ satisfies the measure turnpike property at some steady state. Given any $s\in \sigma^+(A)\cup\sigma^0(A)$ and 
$$
x\in \ker(sI-A)\cap\ker C,
$$
we have that $C\calT_tx=e^{st}Cx=0$ for any $t\in [0,\infty)$. So, by definition, $x\in U^{\infty}$.
Since Theorem~\ref{necess1}~(a) implies that $\calT$ is exponentially stable on $U^{\infty}$, we deduce that $T_tx\to0$ as $t\to \infty$. On the other hand, since $\Re s\geq 0$, 
\begin{equation*}
\|T_tx\|=|e^{st}|\|x\|\geq \|x\|,\;\;\forall t\in[0,\infty).
\end{equation*}
Hence we may conclude that $x=0$. Thanks to Lemma~\ref{Hautus}, we deduce that the pair $(A,C)$ is exponentially detectable. Moreover, Theorem~\ref{necess1}(c) ensures that $(A,B)$ is exponentially stabilizable. Therefore, $(b)\Sra(c)$ holds true.
\end{proof}

\subsection{Proof of Theorem \ref{suffnece3}}
We start with the following technical lemma.
\begin{lem}\label{strictposi}
Suppose problem $(LQ)_T$ satisfies the exponential turnpike property at $(0,0)$. Let $V=\ker[A\;B]$, then the operator $\calP:V\to V$ defined by 
$$
\calP:=\mathbb{P}_{V} \bmat{C^{*}C}{0}{0}{K^{*}K}\bigg|_{V}
$$
is strictly positive.
\end{lem}
\begin{proof}
Let $V$ and $\calP$ be defined as above. Obviously $\calP$ is non-negative. If $\calP$ is not strictly positive, then there exists a sequence $((x_n,u_n))_{n\in\N}\subset V$ such that $\|x_n\|^2+\|u_n\|^2=1$ and
\begin{equation*}
\<\calP\bvec{x_n}{u_n},\bvec{x_n}{u_n}\>=\|Cx_n\|^2+\|Ku_n\|^2\to 0.
\end{equation*}
Since $K$ is coercive, we must have that $u_n\to0$ and $\|x_n\|\to1$ as $n\to\infty$. Since problem $(LQ)_T$ satisfies the exponential turnpike property at $(0,0)$, there exist $M_{\calN},k>0$ such that
\begin{equation}\label{LQexpturn}
\|x^*_T(t,x_0)\|+\|u^*_T(t,x_0)\|\leq M_{\calN}(e^{-kt}+e^{-k(T-t)}),\;\forall t\in[0,T]
\end{equation}
holds for any $\|x_0\|\leq 1$.
Now, fix some sufficiently large $T>0$ such that $M_{\calN}e^{-k\frac{T}{2}}\leq \frac{1}{4}$. Notice that since $(x_n,u_n)$ is a steady state,
\begin{equation}\label{posiest1}
x_n=\calT_tx_n+\Phi_tu_n,\;\forall n\in\N,\;t\in[0,T].
\end{equation}
Since $u_n\to0$, we have that $\Phi_t u_n\to0$ as $n\to\infty$ for any $t\in[0,T]$. From~\eqref{posiest1}, we obtain that $x_n-\calT_tx_n\to 0$ as $n\to\infty$ for any $t\in[0,T]$. Recall from our selection of $(x_n,u_n)$ that $\|x_n\|\to1$, so $
\|\calT_tx_n\|\to1$ as $n\to\infty$ for any $t\in[0,T]$.
Besides, notice that
\begin{equation*}
\left\|x^*_T\left(t,x_n\right)\right\|=\left\|\calT_tx_n+\Phi_tu^*_T\left(\cdot,x_n\right)\right\|\geq \|\calT_tx_n\|-\left\|\Phi_tu^*_T\left(\cdot,x_n\right)\right\|,\;\forall t\in[0,T].
\end{equation*}
Now take $t=\frac{T}{2}$. Combining the above estimate with \eqref{LQexpturn}, we obtain that
\begin{align*}
\liminf_{n\to\infty}M\left\|u^*_T\left(\cdot,x_n\right)\right\|_{L^2}&\geq
\liminf_{n\to\infty}\left\|\Phi_\frac{T}{2}u^*_T\left(\cdot,x_n\right)\right\|\\
&\geq\liminf_{n\to\infty}\bigg(\bigg\|\calT_\frac{T}{2}x_n\bigg\|-\left\|x^*_T\left(\frac{T}{2},x_n\right)\right\|\bigg)\\
&\geq 1-2M_{\calN}e^{-k\frac{T}{2}}\geq \frac{1}{2},
\end{align*}
where $M>0$ in the term on the left-hand side is some constant dependent on $T$. We have thus proved that
\begin{equation*}
\liminf_{n\to\infty}\|u^*_T(\cdot,x_n)\|_{L^2}\geq \frac{1}{2M}.
\end{equation*}
Since $K$ is coercive, there exists some constant $m>0$ such that
\begin{align}\label{posiest2}
\begin{split}
\liminf_{n\to\infty}J_T(x_n,u^*_T(\cdot,x_n)):&=\liminf_{n\to\infty}\int_0^T\|Cx^*_T(t,x_n)\|^2+\|Ku^*_T(t,x_n)\|^2dt\\
&\geq \liminf_{n\to\infty}m\|u^*_T(\cdot,x_n)\|^2_{L^2}\geq \frac{m}{4M^2}.
\end{split}
\end{align}
On the other hand, considering the constant input $u\equiv u_n$ and trajectory $x\equiv x_n$, we obtain
\begin{equation*}
J_T(x_n,u^*_T(\cdot,x_n))\leq J_T(x_n,u)=T(\|Cx_n\|^2+\|Ku_n\|^2)\to 0
\end{equation*}
as $n\to\infty$. This gives a contradiction to \eqref{posiest2}, thus proving our claim.
\end{proof}
We are now well prepared to prove Theorem \ref{suffnece3}.
\begin{proof}[Proof of Theorem \ref{suffnece3}]
Suppose problem $(GLQ)_T$ satisfies the exponential turnpike property at some steady state $(x_e,u_e)$ and $\calN$ is some bounded neighborhood of $0$ in $\calH$. 

By the exponential turnpike property, there exists some positive constants $M_1$ and $k>0$ such that for any $y\in2\calN$ and $t\in[0,T]$,
\begin{equation}\label{iffexpest}
\|x^*_T(t,x_e+y)-x_e\|+\|u^*_T(t,x_e+y)-u_e\|\leq M_{1}(e^{-kt}+e^{-k(T-t)}).
\end{equation}
Observe that for any $T>0$ and $x_0\in\calH$, the trajectory
$$
\wt x^*_T(\cdot,x_0):=x^*_T(\cdot,x_e+2x_0)-x^*_T(\cdot,x_e+x_0)
$$
defined on $[0,T]$ satisfies $\wt x^*_T(0,x_0)=x_0$ and
\begin{align*}
\frac{d\wt x^*_T(t,x_0)}{dt}=&A\wt x^*_T(t,x_0)+B(u^*_T(t,x_e+2x_0)-u^*_T(t,x_e+x_0))
%\\
%=&(A-B(K^*K)^{-1}B^*P(T-t))\wt x^*_T(t,x_0)
\end{align*}
for any $t\in[0,T]$. Similarly, we define
\begin{equation*}
\wt u^*_T(\cdot,x_0):=u^*_T(\cdot,x_e+2x_0)-u^*_T(\cdot,x_e+x_0)
\end{equation*}
on $[0,T]$. Thanks to Corollary~\ref{necess2}, for any $t\in[0,T]$,
\begin{equation*}
u^*_T(t,x_e+2x_0)-u^*_T(t,x_e+x_0)=-(K^*K)^{-1}B^*P(T-t)\wt x^*_T(t,x_0).
\end{equation*}
Therefore, we deduce that
\begin{align*}
\frac{d\wt x^*_T(t,x_0)}{dt}=&(A-B(K^*K)^{-1}B^*P(T-t))\wt x^*_T(t,x_0)\, ,
\end{align*}
which implies that $\wt x^*_T(\cdot,x_0)$ is the optimal trajectory of problem $(LQ)_T$ corresponding to time horizon $T$ and initial condition $x_0$, with corresponding optimal control given by $\wt u^*_T(\cdot,x_0)$.

We now notice that, for any $t\in[0,T]$,
\begin{align*}
\|\wt x^*_T(t,x_0)\|&=\|x^*_T(t,x_e+2x_0)-x^*_T(t,x_e+x_0)\|\\
&\leq \|x^*_T(t,x_e+2x_0)-x_e\|+\|x^*_T(t,x_e+x_0)-x_e\|,
\end{align*}
and
\begin{align*}
\|\wt u^*_T(t,x_0)\|&=\|u^*_T(t,x_e+2x_0)-u^*_T(t,x_e+x_0)\|\\
&\leq \|u^*_T(t,x_e+2x_0)-u_e\|+\|u^*_T(t,x_e+x_0)-u_e\|.
\end{align*}
Combining these two estimates with equation \eqref{iffexpest}, we obtain that
\begin{equation*}
\|\wt x^*_T(t,x_0)\|+\|\wt u^*_T(t,x_0)\|\leq 2M_{1}(e^{-kt}+e^{-k(T-t)}),\;\;\forall t\in[0,T].
\end{equation*}
Since $\calN,T$ and $x_0$ can all be arbitrarily chosen, the exponential turnpike property is satisfied for problem $(LQ)_T$ at $(0,0)$.

Conversely, assume that problem $(LQ)_T$ satisfies the exponential turnpike property at $(0,0)$. Recall that the optimal trajectory $x^*$ of problem $(LQ)_T$ corresponding to time horizon $T>0$ and initial condition $x_0\in\calH$ is given by
\begin{equation*}
x^*(t)=U_T(t,0)x_0,\;\;\forall t\in[0,T].
\end{equation*}
Let $\calN$ be the closed unit ball in $\calH$ centered at $0$. Then the exponential turnpike property implies that there exist some positive constants $M_0$ and $k_0$ such that
\begin{equation*}
\|U_T(t,0)x_0\|\leq M_0(e^{-k_0t}+e^{-k_0(T-t)}),\;\;\forall x_0\in\calN,T>0,t\in[0,T].
\end{equation*}
Equivalently, 
\begin{equation}\label{iffest1}
\|U_T(t,0)\|\leq M_0(e^{-k_0t}+e^{-k_0(T-t)}),\;\;\forall T>0,t\in[0,T].
\end{equation}
Besides, notice that, for any $T>0$
\begin{equation*}
U_{T}(T,0)=U_{T}\left(T,\frac{T}{2}\right)U_{T}\left(\frac{T}{2},0\right)=U_{\frac{T}{2}}\left(\frac{T}{2},0\right)U_{T}\left(\frac{T}{2},0\right),
\end{equation*}
so
\begin{equation}\label{iffest2}
\|U_{T}(T,0)\|\leq M_0(e^{-k_0\frac{T}{2}}+1)2M_0e^{-k_0\frac{T}{2}}\leq 4M_0^2e^{-\frac{k_0}{2}T}.
\end{equation}
Let $V=\ker[A\;B]$. Since Lemma~\ref{strictposi} implies
\begin{equation*}
\calP:=\mathbb{P}_{V} \bmat{C^{*}C}{0}{0}{K^{*}K}\bigg|_{V}
\end{equation*}
is strictly positive, $\calP$ is invertible (see, e.g., \cite[Remark 3.3.4]{TusWei}). Then $\ran\calP=V$ and there exists some $(x_e,u_e)\in V$ such that
\begin{equation*}
\calP\bvec{x_e}{u_e}=\mathbb{P}_{V} \bmat{C^{*}C}{0}{0}{K^{*}K}\bvec{x_e}{u_e}=-\mathbb{P}_{V}\bvec{z}{v}\in V,
\end{equation*}
which further implies that
\begin{equation*}
\bvec{C^{*}Cx_e+z}{K^{*}Ku_e+v}\in V^\perp=\ol{\ran}\bvec{A^*}{B^*}.
\end{equation*}
Then it follows from Theorem \ref{necess1}(c) and Lemma \ref{test} that $\ker\bvec{A^*}{B^*}=\{0\}$ and $\ran\bvec{A^*}{B^*}$ is closed, so there exists a unique $w\in D(A^*)$ satisfying
\begin{equation*}
\bvec{A^*}{B^*}w=\bvec{z+C^*Cx_e}{v+K^*Ku_e}.
\end{equation*}
Following line by line the
proof of Lemma \ref{explicitoc}, we deduce that the optimal control of problem $(GLQ)_T$ is given in a feedback law form by
\begin{align*}
\begin{split}
u_T^*(t,x_0)-u_e=-(K^*&K)^{-1}B^*P(T-t)(x_T^*(t,x_0)-x_e)\\
&-(K^*K)^{-1}B^*(U_{T-t}(T-t,0))^*w,\;\;\forall t\in [0,T].
\end{split}
\end{align*}
It then follows that, for any $t\in[0,T]$,
\begin{align*}
x_T^*(t,x_0)-x_e=&\ \ U_T(t,0)(x_0-x_e)\\
&\ -\int_0^tU_{T-s}(t-s,0)B(K^*K)^{-1}B^*(U_{T-s}(T-s,0))^*w\;ds.
\end{align*}
Simple considerations about \eqref{iffest1}, \eqref{iffest2} and the above equation show that there exists some positive constant $M_1$ such that, for any $t\in[0,T]$,
\begin{align*}
\|x_T^*(t,x_0)-x_e\|\leq M_0&(e^{-k_0t}+e^{-k_0(T-t)})\|x_0-x_e\|\\
&+\int_0^tM_1 (e^{-k_0(t-s)}+e^{-k_0(T-t)})e^{-\frac{k_0}{2}(T-s)}ds.
\end{align*}
A straightforward computation shows that
\begin{equation}\label{iffexpx}
\|x_T^*(t,x_0)-x_e\|\leq M_0 (e^{-k_0t}+e^{-k_0(T-t)})\|x_0-x_e\|+M_2e^{-\frac{k_0}{2}(T-t)},\;\;\forall t\in[0,T]
\end{equation}
holds with some constant $M_2>0$ (independent of $T$).
Finally, recall that
\begin{align*}
\begin{split}
u_T^*(t,x_0)-u_e=-(K^*&K)^{-1}B^*P(T-t)(x_T^*(t,x_0)-x_e)\\
&-(K^*K)^{-1}B^*(U_{T-t}(T-t,0))^*w,\;\;\forall t\in [0,T].
\end{split}
\end{align*}
By the uniform boundedness principle, \eqref{iffest2} and \eqref{iffexpx}, there exist some $M_3,M_4>0$ (independent of $T$) such that
\begin{equation}\label{iffexpu}
\|u_T^*(t,x_0)-u_e\|\leq M_3 (e^{-k_0t}+e^{-k_0(T-t)})\|x_0-x_e\|+M_4e^{-\frac{k_0}{2}(T-t)},\;\;\forall t\in[0,T].
\end{equation}

The combination of \eqref{iffexpx} and \eqref{iffexpu} ensures that problem $(GLQ)_T$ satisfies the exponential turnpike property at $(x_e,u_e)$.
\end{proof}

\section*{Acknowledgements}
The authors thank the referees for their valuable comments, that helped improving the first version of the paper.

The authors acknowledge the support of the Natural Sciences and Engineering Research Council of Canada (NSERC), funding reference number RGPIN-2021-02632.

\bibliographystyle{abbrv}
\bibliography{references.bib}

\appendix
\section{The closed form solution of the optimal control}

\begin{lem}\label{explicitoc} For any $T>0$ and $x_0\in\calH$, if there exists a vector $w\in D(A^*)$ and a steady state $(x_e,u_e)$ such that
\begin{equation}\label{eq:app1}
\bvec{A^*}{B^*}w=\bvec{z+C^*Cx_e}{v+K^*Ku_e},
\end{equation}
then the optimal pair $(x^*_T(\cdot,x_0),u^*_T(\cdot,x_0))$ of problem $(GLQ)_T$ satisfies
\begin{align}\label{ocfeedback}
\begin{split}
u_T^*(t,x_0)-u_e=-(K^*&K)^{-1}B^*P(T-t)(x_T^*(t,x_0)-x_e)\\
&-(K^*K)^{-1}B^*(U_{T-t}(T-t,0))^*w,\;\;\forall t\in [0,T].
\end{split}
\end{align}
\end{lem}
\begin{proof}
For any $T>0$, we claim that $(U_T(T,0))^*$ is the evolution operator at time $T$ of the following evolution system:
\begin{equation*}
\addtolength{\jot}{5pt}
\left\{
\begin{alignedat}{2}
&\frac{dp(t)}{dt}=(A^*-P(t)B(K^*K)^{-1}B^*)p(t),\\
&p(0)=p_0\in\calH.
\end{alignedat}
\right.
\end{equation*}
In other words, we have that $(U_T(T,0))^*p_0=p(T)$.

The proof is based on Yosida approximations. Let $A_n:=nA(nI-A)^{-1}\in \calL(\calH)$ denote the Yosida approximation of $A$ for sufficiently large $n\in\N$. Assume that $x_0,p_0\in \calH$. For each $n$, we let $x_n$ denote the solution of problem
\begin{equation*}
\addtolength{\jot}{5pt}
\left\{
\begin{alignedat}{2}
&\frac{dx_n(t)}{dt}=(A_n-B(K^*K)^{-1}B^*P(T-t))x_n(t),\;\;t\in[0,T]\\
&x_n(0)=x_0,
\end{alignedat}
\right.
\end{equation*}
and $p_n$ denote the solution of problem
\begin{equation*}
\addtolength{\jot}{5pt}
\left\{
\begin{alignedat}{2}
&\frac{dp_n(t)}{dt}=(A_n^*-P(t)B(K^*K)^{-1}B^*)p_n(t),\;\;t\in[0,T]\\
&p_n(0)=p_0.
\end{alignedat}
\right.
\end{equation*}
Since $A_n$ is bounded, we can easily verify (by showing that the derivative is $0$) that
\begin{equation*}
\<x_0,p_n(T)\>=\<x_n(T),p_0\>.
\end{equation*}
By \cite[Part II, Chapter 1, Proposition 3.4]{Bensou}, $p_n(T)\to p(T)$ and $x_n(T)\to U_T(T,0)x_0$ as $n\to\infty$. This further implies that
\begin{equation*}
\<x_0,p(T)\>=\<U_T(T,0)x_0,p_0\>=\<x_0,(U_T(T,0))^*p_0\>.
\end{equation*}
Since $x_0,p_0\in\calH$ can be chosen arbitrarily, our claim then follows.

Now fix some $T>0$ and define $p(\cdot)=(U_{\cdot}(\cdot,0))^*w$ on $[0,T]$. Since $w\in D(A^*)$, by \cite[Part II, Chapter 1, Proposition 3.5]{Bensou} we have that
$$
p(\cdot)\in C([0,T],D(A^*))\cap C^1([0,T],\calH).
$$
We then deduce that
\begin{align*}
\frac{dp(T-t)}{dt}=-(A^*-P(T-t)B(K^*K)^{-1}B^*)p(T-t)\;\;\text{in}\;\;\calH,\;\;\forall t\in[0,T].
\end{align*}
Let $x$ be the solution of problem \eqref{sy} corresponding to input $u\in L^2(0,T;\calU)$ and initial condition $x_0\in\calH$. Then, thanks to assumption~\eqref{eq:app1},
\begin{align*}
\int_0^T\<z,x(t)&\>+\<v,u(t)\>dt\\
&=\int_0^T\<w,Ax(t)+Bu(t)\>_{\calH_1^d,\calH_{-1}}-\<Cx_e,Cx(t)\>-\<Ku_e,Ku(t)\>dt\\
&=\<w,x(T)-x_0\>+\int_0^T-\<Cx_e,Cx(t)\>-\<Ku_e,Ku\>dt.
\end{align*}
Combining this with \eqref{Ptransform} shows that
\begin{align}\label{covGLQ}
\begin{split}
J_T(x_0,u)=\int_0^T&\|K\{(u(t)-u_e)+(K^*K)^{-1}B^*P(T-t)(x(t)-x_e)\}\|^2dt\\
&-J_T(x_e,u_e)+\<P(T)(x_0-x_e),(x_0-x_e)\>+2\Re\<w,x(T)-x_0\>. 
\end{split}
\end{align}
On the other hand, we have
\begin{align*}
\begin{split}
2\Re\<&x(T),w\>-2\Re\left\<x_0,p(T)\right\>\\
&=\Re\int_0^T2\left\<Ax(t)+Bu(t),p(T-t)\right\>_{H_{-1},H_{1}^d}\\
&\qquad\qquad+2\left\<x(t),-(A^*-P(T-t)B(K^*K)^{-1}B^*)p(T-t)\right\>dt\\
&=\Re\int_0^T2\left\<u(t),B^*p(T-t)\right\>\\
&\qquad\qquad+2\left\<(K^*K)^{-1}B^*P(T-t)x(t),B^*p(T-t)\right\>dt.
\end{split}
\end{align*}
This, together with \eqref{covGLQ} implies that
\begin{align*}
J_T(x_0,u)&=\int_0^T\|K\{(u(t)-u_e)+(K^*K)^{-1}B^*P(T-t)(x(t)-x_e)\}\|^2dt\\
&\qquad\qquad+2\Re\<w,x(T)\>+M_0\\
&=\int_0^{T} \|K\{(u(t)-u_e)+(K^*K)^{-1}B^*P(T-t)(x(t)-x_e)\}\|^2\\
&\qquad\qquad+2\Re\left\<(K^*K)^{-1}B^*P(T-t)(x(t)-x_e),B^*p(T-t)\right\>\\
&\qquad\qquad+2\Re\left\<(u(t)-u_e),B^*p(T-t)\right\>dt+M_1\\
&=\int_0^{T}\left\|K\left\{(u(t)-u_e)+(K^*K)^{-1}B^*\left[P(T-t)(x(t)-x_e)+p(T-t)\right]\right\}\right\|^2dt\\
&\qquad\qquad+M_2
\end{align*}
where $M_0$, $M_1$ and $M_2\in \R$ are constants independent of $u$.

This implies that, if the closed loop system corresponding to the feedback law
\begin{align*}
u(t)-u_e=-(K^*&K)^{-1}B^*P(T-t)(x(t)-x_e)\\
&-(K^*K)^{-1}B^*(U_{T-t}(T-t,0))^*w,\;\;\forall t\in [0,T],
\end{align*}
admits a solution in $L^2(0,T;\calH)$, then this solution is optimal. 

By \cite[Part II, Chapter 1, Proposition 3.5]{Bensou}, this problem does admit a unique solution $x$ in $L^2(0,T;\calH)$. So, the optimal pair verifies equation \eqref{ocfeedback}.
\end{proof}

\end{document}